\documentclass[11pt]{amsart}
\usepackage{hyperref}
\usepackage{amsmath, amssymb, amsthm, bm, mathtools, enumitem, caption}
\usepackage[noabbrev,capitalize]{cleveref}
\usepackage{adjustbox}
\crefname{equation}{}{}

\usepackage{graphicx, tikz}
\usepackage[textheight=8.7in, textwidth=6.7in]{geometry}

\usepackage{microtype}

\newtheorem{theorem}{Theorem}[section]

\newtheorem*{conjecture*}{Conjecture}

\theoremstyle{definition}

\theoremstyle{remark}
\newtheorem*{remark}{Remark}

\DeclarePairedDelimiter\abs{\lvert}{\rvert}

\newcommand{\p}[1]{\partial #1}

\newcommand{\QQ}{\mathbb Q}

\newcommand{\ZZ}{\mathbb Z}
\renewcommand{\abs}[1]{\left\vert #1 \right \vert}

\newtheorem*{conjectureeh}{Conjecture EH}
\title[]{Euler--Kronecker constants for cyclotomic fields}
\date{\today}
\thanks{2010 {\it{Mathematics Subject Classification.}} 11M06, 11N37, 11R18, 11R42, 11Y60}
\keywords{Cyclotomic fields; Elliott-Halberstam Conjecture; Euler-Kronecker constants}
\author{Letong Hong \and Ken Ono \and Shengtong Zhang}
\address{Dept. of Mathematics, Massachusetts Institute of Technology, Cambridge, MA 02139}
\email{clhong@mit.edu}
\address{Dept. of Mathematics, University of Virginia, Charlottesville, VA 22904}
\email{ko5wk@virginia.edu}

\address{Dept. of Mathematics, Massachusetts Institute of Technology, Cambridge, MA 02139}
\email{stzh1555@mit.edu}

\begin{document}
\begin{abstract}  The Euler-Mascheroni constant
$\gamma=0.5772\dots\!$ is the $K=\QQ$ example of an Euler-Kronecker constant $\gamma_K$ of a number field $K.$ In this note we consider the size of the $\gamma_q=\gamma_{K_q}$ for cyclotomic fields $K_q:=\QQ(\zeta_q).$
Assuming the Elliott-Halberstam Conjecture (EH), we prove uniformly in $Q$ that
$$\frac{1}{Q}\sum_{Q<q\le 2Q} \abs{\gamma_q - \log q} = o(\log Q).$$
In other words,  under EH the  $\gamma_q / \log q$ in these ranges converge to the one point distribution at $1$. This theorem refines and extends a previous result of Ford, Luca, and Moree for prime $q.$  The proof of this result is a straightforward modification of earlier work of Fouvry under the assumption of EH.

\end{abstract}

\maketitle
\section{Introduction}

For a number field $K$, the {\it Euler-Kronecker constant} $\gamma_K$ is given by
$$\gamma_K:=\lim_{s \to 1^+}\left(\frac{\zeta_K'(s)}{\zeta_K(s)} + \frac{1}{s-1}\right),$$ where $\zeta_K(s)$ is  the Dedekind zeta-function for $K.$
The  Euler-Mascheroni constant
$\gamma=0.5772\dots\!$ is the $K=\QQ$ case, where $\zeta_{\QQ}(s)=\zeta(s)$ is the Riemann zeta-function.
We consider the constants $\gamma_q=\gamma_{K_q}$ for cyclotomic fields $K_q:=\QQ(\zeta_q)$, where $q\in \ZZ^{+}$ and $\zeta_q$ is a primitive $q$th root of unity.

The recent interest in the distribution of the $\gamma_q$ is inspired by work of Ihara \cite{Ihara, Ihara2}. He proposed,
 for every $\varepsilon > 0,$ that there is a $Q (\varepsilon)$ for which
 $$ (c_1-\varepsilon)\log q\le \gamma_q\le (c_2+\varepsilon)\log q$$
 for every integer $q \ge Q (\epsilon)$, where  $0 < c_1 \le c_2 < 2$ are absolute constants.
This conjecture was disproved by Ford, Luca and Moree in \cite{FLM} assuming a strong form of the Hardy--Littlewood $k$-tuple Conjecture.
However, assuming the Elliott-Halberstam Conjecture (see \cite{EH}), these same authors also proved that the conjecture holds for almost all primes $q,$ with $c_1 = c_2 = 1.$ 
We recall the Elliott-Halberstam Conjecture as formulated in terms of the Von Mangoldt function $\Lambda(n),$ the Chebyshev function $\psi(x),$ and Euler's totient function $\varphi(n).$

\begin{conjectureeh} If we let
 $$E(x;m,a):=\sum_{\substack{p\equiv a\!\!\pmod m\\ p\leq x\ {\text {\rm prime}}}}
  \Lambda(p) -\frac{\psi(x)}{\varphi(m)},$$ 
 then for every $\varepsilon>0$ and $A>0$, we have
 $$\sum_{m\le x^{1-\varepsilon}}\max_{(a,m) = 1}\abs{E(x;m,a)}\ll_{A,\varepsilon} \frac{x}{(\log x)^A}.$$
\end{conjectureeh}

Assuming EH, Ford, Luca, and Moree proved (see Theorem 6 (i) of \cite{FLM}),
for every $\varepsilon > 0$, that
 $$1-\varepsilon < \frac{\gamma_q}{\log q} <1+\varepsilon $$
for almost all primes $q$ (that is, the number of exceptional $q \le x$ is $o(\pi(x))$ as $x\to \infty$).
Here we extend and refine this result to all integers $q.$

\begin{theorem}\label{main}
Under EH,  for $Q\rightarrow +\infty$ we have
$$\frac{1}{Q}\sum_{Q<q\le 2Q} \abs{\gamma_q - \log q} = o(\log Q),$$
where the sum is over integers $q$.
\end{theorem}

\begin{remark}
Theorem~\ref{main} shows that EH implies that the distribution of $\gamma_q / \log q$ in $[Q, 2Q]$
 converges to the one point distribution supported on $1$.
\end{remark}

To prove Theorem~\ref{main}, we use of work of Fouvry \cite{Fouvry} that allowed him to unconditionally prove that
 $$\frac{1}{Q}\sum_{Q<q\le 2Q}\gamma_q=\log Q+O(\log\log Q).$$  Our conditional result is a point-wise refinement
 of Fouvry's asymptotic formula under EH.

\section{Acknowledgements}
\noindent
The authors acknowledge the Thomas Jefferson Fund and the NSF (DMS-2002265 and DMS-2055118) for their generous support. The authors thanks Pieter Moree for helpful discussions regarding his work with Ford and Luca. We thank the referee for suggestions that improved the exposition in this paper.

\section{Proof of \cref{main}}

For brevity, we shall assume that the reader is familiar with Fouvry's paper \cite{Fouvry}.
The key formula is (see (3) of \cite{Fouvry}) the following expression for $\gamma_q$ in terms of logarithmic derivatives of Dirichlet $L$-functions:
\begin{equation}\label{formula}
\gamma_q = \gamma + \sum_{1<q^*|q\ }\sum_{\chi^*\bmod q^*}\frac{L'(1,\chi^*)}{L(1,\chi^*)}.
\end{equation}
Here the inner sum runs over the primitive Dirichlet characters $\chi^*$ modulo $q^*$.

We follow the strategy and notation in \cite{Fouvry}, which makes use of the modified Chebyshev function
\begin{displaymath}
\psi(x;q,a):=\sum_{\substack{n\leq x\\ n\equiv a\!\!\!\!\pmod q}} \Lambda(n),
\end{displaymath}
and the integral
$$
\Phi_{\chi^*}(x):=\frac{1}{x-1}\int_{1}^x \left(\sum_{n\leq t}\frac{\Lambda(n)}{n}\chi^*(n)\right) dt.
$$
However, we replace the sums $\Gamma_i(Q)$ and $\Gamma_{1,j}(Q)$ defined in \cite{Fouvry} with the pointwise terms $\gamma_i(q)$ and $\gamma_{1,j}(q)$. Following the approach in \cite{Fouvry},  which is based on (\ref{formula}), we have
$$\gamma_q = \gamma + A(q) + B(q) - \gamma_2(q) - \gamma_3(q) - (\gamma_{1,1}(q) + \gamma_{1,2}(q) + \gamma_{1,3}(q)),$$
where
$$A(q)=\sum_{q^*|q} \sum_{\chi^*\mod q^*} \frac{L'}{L}(1,\chi^*)+\Phi_{\chi^*}(x),$$
$$B(q)=\sum_{\substack{\chi \bmod q\\ \chi\neq \chi_0}} \Phi_{\chi}(x)-\sum_{q^*|q}\sum_{\chi^*\mod q^*} \Phi_{\chi^*}(x),$$
$$\gamma_2(q)=\frac{1}{x-1}\int_1^x \frac{\varphi(q)\psi(t;q,1)-\psi(t)}{t} dt,$$
$$\gamma_3(q)=\frac{1}{x-1}\int_1^x \sum_{\substack{n\leq t\\(n,q) \neq 1}} \frac{\Lambda(n)}{n}dt,$$
$$\gamma_{1,1}(q)=\frac{1}{x-1}\int_1^x\int_1^{\min(q,t)}\left(\frac{\varphi(q)\psi(u;q,1)-\psi(u)}{u^2}du\right)dt,$$
$$\gamma_{1,2}(q)=\frac{1}{x-1}\int_1^x\int_{\min(q,t)}^{\min(x_1,t)}\left(\frac{\varphi(q)\psi(u;q,1)-\psi(u)}{u^2}du\right)dt,$$
$$\gamma_{1,3}(q)=\frac{1}{x-1}\int_1^x\int_{\min(x_1,t)}^t \left(\frac{\varphi(q)\psi(u;q,1)-\psi(u)}{u^2}du\right)dt.$$
To complete the proof, for $\varepsilon >0$ we let
  $x := q^{100}$ and  $x_1 := q^{1 + \varepsilon}$. Apart from $\gamma_{1,1}(q),$ which gives the  $-\log q$ terms in Theorem~\ref{main}, we shall show that these summands  are all small.

\medskip
\noindent
\textbf{Estimation of $A(q)$}: By Proposition 1 and Remark (i) of \cite{Fouvry}, we have
$$\sum_{q=Q}^{2Q} \abs{A(q)}=O(Q).$$
\textbf{Estimation of $B(q)$}: For $B(q)$, by equation (26) and Lemma 3 of \cite{Fouvry},  we simplify 
\begin{align*}
B(q)&=-\frac{1}{x-1}\int_{1}^x \sum_{q^*|q}\ \sum_{\chi^*\!\! \!\mod q^*} \sum_{\substack{n \leq t\\(n,q)>1}} \frac{\Lambda(n)\chi^*(n)}{n} dt \\    
&= -\frac{1}{x-1}\int_{1}^x \sum_{q^*|q}\ \sum_{\chi^*\!\! \!\mod q^*}\ \sum_{\substack{p^v\leq t\\p|q}} \frac{\log p\cdot\chi^*(p^v)}{p^v} dt \\
&= -\frac{1}{x-1}\int_1^x \sum_{q^*|q}\ \sum_{\substack{p^v\leq t \\ p|q\\p\nmid q^*}}\  \sum_{\substack{d|(p^v-1,q^*)}}\frac{\log p}{p^v} \cdot \varphi(d) \mu\left(\frac{q^*}{d}\right) dt \\
&= -\frac{1}{x-1}\int_1^x  \sum_{\substack{p^v\leq t\\p|q}}\ \sum_{\substack{d|p^v-1}}\frac{\log p}{p^v}\cdot\varphi(d) \sum_{\substack{q^*|q\\d|q^*\\p\nmid q^*}}\mu\left(\frac{q^*}{d}\right) dt.
\end{align*}
We note that the innermost sum
$$\sum_{\substack{q^*|q\\d|q^*\\p\nmid q^*}}\mu\left(\frac{q^*}{d}\right)$$
is always $0$ or $1$, so we conclude that $B(q) \leq 0$ for any $q$. Proposition 2 of \cite{Fouvry} gives
$$\sum_{q=Q}^{2Q}B(q)=O(Q),$$
and so we have
$$\sum_{q=Q}^{2Q}\abs{B(q)}=O(Q).$$
\textbf{Estimation of $\gamma_2(q)$}: By Lemma 8 of \cite{Fouvry}, uniformly in $Q$ with $u\geq 1,$ we have
$$\sum_{q=Q}^{2Q} \psi(u;q,1) \ll u.$$
Therefore, we have that
$$\sum_{q=Q}^{2Q} \abs{\varphi(q)\psi(t;q,1)-\psi(t)}=O(Qt),$$
and so we conclude that
$$\sum_{q=Q}^{2Q} \abs{\gamma_{2}(q)} = O(Q).$$
\textbf{Estimation of $\gamma_3(q)$}: By definition, $\gamma_3$ is positive, so by equation (36) of \cite{Fouvry}, we have
$$\sum_{q=Q}^{2Q} \abs{\gamma_{3}(q)}=O(Q).$$
\textbf{Estimation of $\gamma_{1,1}(q)$}:
 Since $\psi(u;q,1)=0$ for $u<q$, we have $$\gamma_{1,1}(q)=-\frac{1}{x-1}\int_1^x \left(\int_1^{\min(q,t)}\frac{\psi(u)}{u^2} du \right) dt.$$ Dividing both sides of Equation (41) of \cite{Fouvry} by $Q$, we have
$$\gamma_{1,1}(q)=-\log q + O(1).$$
\textbf{Estimation of $\gamma_{1,2}(q)$}: By the same proof as  equation (42) of \cite{Fouvry}, we have
$$\sum_{q=Q}^{2Q} \abs{\gamma_{1,2}(q)} \ll \varepsilon Q\log Q.$$
Summing the above estimates, we conclude unconditionally that
$$\frac{1}{Q}\sum_{q=Q}^{2Q} \abs{\gamma_q-\log q}=\frac{1}{Q}\sum_{q=Q}^{2Q} \abs{\gamma_{1,3}(q)} + O(\varepsilon \log Q).$$
\textbf{Estimation of $\gamma_{1,3}(q)$}: If we assume Conjecture EH holds, then we have (as in Lemma 7 of \cite{Fouvry}) that
$$\sum_{\substack{q\le 2Q\\(q,a)=1}}\varphi(q)\abs{\psi(x;q,a)-\frac{\psi(x)}{\varphi(q)}}=O_A\big(Qx(\log x)^{-A+2}\big).$$
Therefore, we find that
$$\frac{1}{Q} \sum_{q=Q}^{2Q}\abs{\gamma_{1,3}(q)}=O_{\epsilon,A}(\log^{-A} Q).$$ 
By combining these estimates, we obtain the main result
$$\frac{1}{Q}\sum_{q=Q}^{2Q}\abs{\gamma_q-\log q}=o(\log Q),$$

\end{document}